\documentclass{amsart}
\usepackage{amssymb}
\newtheorem{Theo}{Theorem}

\newtheorem{Lem}[Theo]{Lemma}

\newcommand{\Q}{\mathbb{Q}}

\newcommand{\F}{\mathbb{F}}
\newcommand{\jj}{\mathfrak{j}}
\newcommand{\qq}{\mathfrak{q}}
\normalsize
\begin{document}
\title{On a criterion for Catalan's Conjecture}
\author{Jan-Christoph Schlage-Puchta}
\begin{abstract}
We give a new proof of a theorem of P. Mih\v ailescu which states that 
the equation $x^p-y^q=1$ is unsolvable with $x, y$ integral and $p, q$ odd
primes, unless the congruences $p^q\equiv p\pmod{q^2}$ and
$q^p\equiv q\pmod{p^2}$. hold.
\end{abstract}
\maketitle
MSC-index: 11D61, 11R18
Improving criterions for Catalan's equation by Inkeri\cite{Ink},
Mignotte\cite{Mig}, Schwarz\cite{Sch} and Steiner\cite{Ste},
Mihailescu\cite{Mih} proved the following theorem.

\begin{Theo}
Let $p, q$ be odd prime numbers. Assume that $p^q\not\equiv p\pmod{q^2}$ or
$q^p\not\equiv q\pmod{p^2}$. Then the equation $x^p-y^q=1$ has no nontrivial
integer solutions.
\end{Theo}

Here we will give a different proof of this theorem. More precisely, we will
show the following statement.

\begin{Theo}
Let $p, q$ be odd prime numbers, and assume that the equation $x^p-y^q=1$ has
some nontrivial solution. Then we have either $q^2|p^q-p$ or the $q$-rank of
the relative class group of the $p$-th cyclotomic field is at least $(p-5)/2$.
\end{Theo}

Note that different from Mihailescu's proof of Theorem 1, we have to make use
of estimates for the relative size of $p$ and $q$ obtained using bounds for
linear forms in logarithms, thus the passage from Theorem 2 to Theorem 1 is by
no means elementary. However, the proof of Theorem 2 makes much less use of
special properties of cyclotomic fields than Mihailescu's proof of Theorem 1,
thus it might be easier to adapt to different situations. 

To deduce Theorem 1 from Theorem 3, it suffices to show that the second
alternative is
impossible. Assume that $x^p-y^q=1$, and that the $q$-rank of the relative
class group of the $p$-th cyclotomic field is at least $(p-5)/2$. This implies
$q^{(p-5)/2}\leq h^-(p)$. The class number $h^-(p)$ was estimated by Masley
and Montgomery\cite{MM}, they showed that for $p>200$ we have
$h^-(p)<(2\pi)^{-p/2} p^{(p+31)/4}$. Thus we get $q<\sqrt{p}$. On the other
hand, Mignotte and Roy\cite{MR} proved, that for $q\geq 3000$ we have
$p\leq 2.77q\log q
(\log p - \log\log q + 2.33)^2$, combining these inequalities and observing
that Mignotte and Roy\cite{MR2} have shown that $q>10^5$, thus
$\log\log q>2.33$, we get $p\leq 1.92\log^6 p$, which implies
$p<6.6\cdot 10^7$, thus $q<\sqrt{p}<8200$ contradicting the lower bound
$q>10^5$ mentioned above.

To prove theorem 3, we follow the lines of \cite{Sch}, incorporating an idea
of Eichler\cite{Eic}. $K$ be the $p$-th
cyclotomic field, $\zeta$ a $p$-th root of unity, $I_K$ the group of fractional
ideals in $K$, $i:K^*\to I_K$ the canonical map $x\mapsto (x)$,
$K^+= \Q(\zeta+\zeta^{-1})$ be the maximal real subfield of $K$, $\mathcal{O}_K$
be the ring of integers of $K$. Denote with $r$ the $q$-rank of the relative
class group of $K$. We begin with a Lemma. $\mathcal Q$ be the set of prime ideals
dividing $q$ in $K$. Choose a primitive root $g$ of $p$ and define
$\sigma\in\mbox{Gal}(K|\Q)$ by the relation $\zeta^\sigma=\zeta^g$.

\begin{Lem}
There is a subgroup $I_0$ of $I_K$ with the following properties:
\begin{enumerate}
\item The prime ideals in $\mathcal Q$ do not appear in the factorization of any
ideal in $I_0$
\item $I_K/(i(K^*)I_0)$ has $q$-rank $r$
\item If $\epsilon\in K^*$ with $(\epsilon)\in I_0$, then
$\epsilon/\overline{\epsilon}$ is a root of unity.
\end{enumerate}
\end{Lem}
{\em Proof:} This is Lemma 1 in \cite{Sch}.

Now assume that $x$ and $y$ are nonzero integers with $x^p-y^q=1$. We have
\nolinebreak\cite{Ink}
\[
\left(\frac{x-\zeta}{1-\zeta}\right) = \jj^q
\]
for some integral ideal $\jj$. The ideal classes with $\jj^q=(1)$ generate an
$r$-dimensional vector space over $\F_q$ in $I_K/(i(K^*)I_0)$, hence there are
integers $a_0, \ldots, a_r$, not all divisible by $q$, such that
$\jj^{a_0+a_1\sigma+\ldots+a_r\sigma^r}$ lies in $i(K^*)I_0$. Thus we get
\[
\left(\frac{x-\zeta}{1-\zeta}\right)^{a_0+a_1\sigma+\ldots+a_r\sigma^r} =
\epsilon\alpha^q
\]
with $(\epsilon)\in I_0$ and $\alpha$ is $\qq$-integral for all prime ideals
$\qq$ dividing $q$, since the left hand side is $\qq$-integral, and
$(\epsilon)$ is not divisible by $\qq$ by condition 1 of Lemma 4. We multiply
this equation with
$(-\zeta^{-1}(1-\zeta))^{a_0+a_1\sigma+\ldots+a_r\sigma^r}$ to get
\begin{equation}
\big(1-x\zeta^{-1}\big)^{a_0+a_1\sigma+\ldots+a_r\sigma^r} = \epsilon'\lambda
\alpha^q
\end{equation}
where $\lambda$ divides some power of $p$, and $\epsilon'$ differs from
$\epsilon$ by some power of $\zeta$, especially $(\epsilon)=(\epsilon')$.

By \cite{Cas}, we have $q|x$, thus the left hand side of (1) can be simplified
$\pmod{q^2}$. We get
\begin{equation}
1-x\left(a_0\zeta^{-1} + a_1\zeta^{-\sigma}+\ldots+a_r\zeta^{-\sigma^r}\right)
\equiv \epsilon'\lambda\alpha^q\pmod{q^2}
\end{equation}
The complex conjugate of the right hand side can be written as
$\zeta^k\epsilon'\lambda\overline{\alpha}^q$, since every $p$-th root of unity
is the $q$-th power of some root of unity, this equals
$\epsilon'\lambda\beta^q$ for some $\beta\in K^*$. Thus if we substract the
complex conjugate of (2), we get
\begin{equation}
x\left(a_0\zeta^{-1} + \ldots + a_r\zeta^{-\sigma^r} -
a_0\zeta^ - \ldots - a_r\zeta^{\sigma^r}\right) \equiv
\epsilon'\lambda(\alpha^q - \beta^q) \pmod{q^2}
\end{equation}
The left hand side of (3) is divisible by $q$, since $x$ is divisible by $q$,
and the bracket is integral. However, $(\epsilon')\in I_0$, and by construction
we have $(\epsilon', q)=(1)$, and $\lambda$ divides some power of $p$, thus
we have $(\lambda, q)=(1)$, too. Hence $q|\alpha^q-\beta^q$, and since $q$ is
unramified, this implies $q^2|\alpha^q-\beta^q$. Hence $q^2$ divides the left
hand side of (3). But $x$ is rational, thus either $q^2|x$, or $q$ divides the
bracket. By \cite{Cas}, we have $x\equiv -(p^{q-1}-1)\pmod{q^2}$, hence the
first possibility implies $q^2|p^q-p$. Thus to prove our theorem, it suffices
to show that the second choice is impossible. 

Assume that
\[
a_0\zeta^{-1} + a_1\zeta^{-\sigma} + \ldots + a_r\zeta^{-\sigma^r} -
a_0\zeta - a_1\zeta^\sigma - \ldots - a_r\zeta^{\sigma^r} = 
q\alpha
\]
This can be written as
\[
a_0X^{\overline{-1}} + a_1X^{\overline{-g}} + \ldots + a_rX^{\overline{-g^r}}
-a_0X-a_1X^g - \ldots - a_rX^{\overline{g^r}} = q F(X) + G(X)\Phi(x)
\]
where $F$ and $G$ are polynomials with rational integer coefficients, $\Phi$
is the $p$-th cyclotomic polynomial, and $\overline{a}$ denotes the least
nonnegative residue $\pmod{p}$ of $a$. The left hand side is of degree
$\leq p-1$, and since we may assume that the leading coefficient of $G$ is
prime to $q$, this implies that $G$ is constant. Further on the left hand side
there are at most $2r+2\leq p-3$ nonvanishing coefficients, thus $G=0$. This
implies that all coefficients on the left hand side vanish $\pmod{q}$. But all
the monomials on the left hand side have different exponents, since otherwise
we would have $g^{s_1}\equiv \pm g^{s_2}\pmod{p}$, which would imply that the
order of $g$ is $\leq 2r\leq p-5$, but $g$ was chosen to be primitive. Hence
all $a_i$ vanish $\pmod{q}$, but this contradicts the choice of the $a_i$ at
the very beginning.

Jan-Christoph Puchta\\
Mathematisches Institut\\
Eckerstra\ss e 1\\
79104 Freiburg\\
Germany\\
jcp@arcade.mathematik.uni-freiburg.de

\begin{thebibliography}{20}
\bibitem{Cas} J. W. S. Cassels, {\em On the equation $a^x-b^y=1$}
Proc. Camb. Philos. Soc. 56, 97-103 (1960)
\bibitem{Eic} M. Eichler, {\em Eine Bemerkung zur Fermatschen Vermutung},
Acta Arith. 11, 129-131 (1965)
\bibitem{Ink} K. Inkeri, {\em On Catalan's Conjecture}, J. Number Theory 34,
142-152 (1990)
\bibitem{MM} J. Masley, H. L. Montgomery, {\em Cyclotomic fields with unique
factorization} J. reine angew. Math. 286/287, 248-256 (1976)
\bibitem{Mig} M. Mignotte, {\em A criterion on Catalan's equation} J. Number
Theory 52,  280-283 (1995)
\bibitem{MR} M. Mignotte, Y. Roy {\em Catalan's equation has no new solution
with either exponent less than 10651} Exp. Math. 4,  259-268 (1995)
\bibitem{MR2} M. Mignotte, Y. Roy, {\em Minorations pour l'equation de Catalan}
C. R. Acad. Sci., Paris, Ser. I 324, 377-380 (1997)
\bibitem{Mih} P. Mih\v ailescu, {\em A class number free criterion for
Catalan's conjecture}, manuscript, Z\"urich (1999)
\bibitem{Sch} W. Schwarz, {\em A note on Catalan's equation} Acta Arith. 72,
277-279 (1995).
\bibitem{Ste} R. Steiner, {\em Class number bounds and Catalan's equation}
Math. Comput. 67, 1317-1322 (1998).
\end{thebibliography}
\end{document}